\documentclass[a4paper]{article}

\usepackage[english]{babel}
\usepackage[utf8]{inputenc}
\usepackage{amsmath}
\usepackage{amsfonts}
\usepackage{graphicx}
\usepackage{epstopdf} 

\title{Desynchronization of Random Dynamical System
under Perturbation by an Intrinsic Noise}

\author{Adrian Jarret\thanks{Summary report for an internship in
the Department of Applied Mathematics, University of Washington, Seattle, under the supervision of Professor Hong Qian.}
\\[10pt]
Department of Applied Mathematics\\
University of Washington, Seattle, WA, 98109, USA\\[1pt]
and\\
\'{E}cole CentraleSup\'{e}lec\\
Universit\'{e} Paris-Saclay,
}

\date{June 18, 2018}

\newtheorem{mydef}{Definition}
\newtheorem{prop}{Property}

\begin{document}
\maketitle

\begin{abstract}
In the theory of random dynamical systems (RDS),
individuals with different initial states follow a same 
law of motion that is stochastically changing with
time | called extrinsic noise.  In the present work,
intrinsic noises for each individual are considered as 
a perturbation to an RDS.  This gives rise to 
random Markov systems (RMS) in which the law of motion
is still stochastically changing with time, but individuals
also exhibit statistically independent variations, with each transition 
having a small probability not to follow the law.  As a 
consequence, two individuals in an RMS system go through 
stochastically distributed periods of synchronization and desynchronization, driven by extrinsic and intrinsic noises
respectively.  We show that in-sync time, {\em e.g.}, escaping from
a random attractor, has a symptotic geometric distribution.
\end{abstract}

\section{Introduction}
\label{sec:introduction}

Might they be biological, physical, numerical, or electronic, most of dynamical phenomena in real world are affected by uncertainties, or noises. The source of a noise can either be {\em extrinsic} that
is same for all the individuals or entities in a system, or {\em intrinsic} thus statistically independent for different individuals. In both cases, noises introduce the possibility of undetermined or random transition between two states in the time evolution of an individual.  Markov chains are a widely used mathematical representation for such {\em stochastic dynamics}. Takahashi presents an interesting study on this kind of Markov chains (MC) using random transition matrices \cite{takahashi1969}, which fits with the noise problem as the transition matrix of the MC is changing with time (random non homogeneous MC). It showed in an accessible manner what can be expected from such mathematical theory.

Recently, random dynamical systems (RDS) have been introduced
as a model for individuals in a system with extrinsic noise.
The topic is well known, and some theorems precisely describe the behavior of such systems; see \cite{Arnold1991} for more details.
Under certain conditions, extrinsic noise has been shown to induce ``synchronizations'' among different individuals in such systems: All the individuals end by being in the same state, having simultaneously the same stochastic transitions.  The limiting behavior is sometime called a
``random attractor''.

In the present work, we study small perturbation of an RDS with
intrinsic noise.  With both extrinsic and intrinsic noises,
the problem is intimately related to the subject known as 
Markov chains in random environments \cite{mcinre}, with
possible applications in Markov decision process \cite{MDP1970}. 
We shall focus on the escape from a random 
attractor.  In a nutshell, synchronizations in the systems is 
driven by the extrinsic noise; but the phenomenons of desynchronization,  when two individuals in the same state at a given time no longer follow a same transition, is caused by independent random disturbances that are different for each individual.  Eventually, it leads to successive synchronizations and desynchronizations.

\section{Definitions and background}
\label{sec:definitions}
In this section, we properly define the concepts that we will use later in this paper. We first give a mathematical definition of the random dynamical systems (RDS), then we introduce a perturbation to extend this definition to random Markov systems (RMS). Finally, we study the behavior of such systems, and how we can in practice understand them with some approximations.

\subsection{Random Dynamical Systems (RDS)}
\label{RDS}
The idea of discrete state, discrete time (dsdt for short) RDS comes from the study of Markov Chains (MC) and its generalization. An homogeneous MC on discrete state space is mainly defined by its transition matrix $M$. Yet, we realize that we can decompose $M$ as a sum of deterministic transition matrices $D_i$ (matrix with only $0$ and $1$, and one entry $1$ only on each row), such that $M = \sum_{i\in\mathcal{I}}{q_iD_i}$, where $\left(q_i\right){i\in\mathcal{I}}$ is a discrete probability distribution on $\left(D_i\right){i\in\mathcal{I}}$. During the realization of a random state following a RDS, at each step we draw a matrix $D_i$ according to this distribution, then the random state follows the (deterministic) transition given by the matrix. This is called an i.i.d. dsdt-RDS. It is easy to show that the expected transition matrix is the transition probability matrix for the MC. So, in average, for a single point motion, one should see the same behavior with both the MC and the RDS.

However, if we consider a two-point motion, i.e. two different random states that follow the same i.i.d. RDS, one can see some dynamics that would be extremely unlikely for a Markov process. Indeed, as long as the two systems are not in the same state at the same time, everything behaves like two realizations of the MC. But let us imagine that, at some point, the two systems enter a same state. Since each transition, once the matrix drawn, are deterministic, the two trajectories would stay together forever. This phenomenon is called synchronization, and is one of the main topics regarding to RDS study. For more complete information on RDS, and on the relation between RDS and MC, you are referred to the article \cite{X.-F.Ye2016}.

Let us take an example and consider a transition matrix on a two-states space
\begin{equation}
M = \left( \begin{array}{cc}
0.2 & 0.8\\
	0.6 & 0.4\\
\end{array}\right).
\label{eq:mean_markov}
\end{equation}
We realize that we can write
$$
M = 0.6 \left( \begin{array}{cc}
0 & 1\\
	1 & 0\\
\end{array}\right) + 
0.2 \left( \begin{array}{cc}
0 & 1\\
	0 & 1\\
\end{array}\right) + 
0.2\left( \begin{array}{cc}
1 & 0\\
	0 & 1\\
\end{array}\right),
$$
but, in the mean time, we can also decompose
$$
M = 0.5 \left( \begin{array}{cc}
0 & 1\\
	1 & 0\\
\end{array}\right) + 
0.3 \left( \begin{array}{cc}
0 & 1\\
	0 & 1\\
\end{array}\right) + 
0.1\left( \begin{array}{cc}
1 & 0\\
	0 & 1\\
\end{array}\right) + 
0.1\left( \begin{array}{cc}
1 & 0\\
	1 & 0\\
\end{array}\right).
$$
So, one can notice that, for each RDS, there exist a unique MC that its mean representation, but a given MC may be represented by different RDS, which potentially different properties (for instance some will synchronize and some not). Again, read \cite{X.-F.Ye2016} for more explanation on the topic. For this reason, in this paper, we only focus on the study of a given RDS, which means a given set $\left(q_i,D_i\right)_{i\in\mathcal{I}}$, and not on any MC. Let us now introduce some notations to define RDS in a more regressed way.

We consider a set of states $\mathcal{S} = \left\{1, 2, \dots, s\right\}$, which size is $s = |\mathcal{S}|$, and the set $\mathcal{D} = \left\{D_1, D_2, \dots \right\}$ of all the deterministic transition matrix of size $s\times s$. We have $|\mathcal{D}| = s^s$. So, each step in the realization of a RDS has its value in
$$
\hat{\Omega} = \mathcal{S}\times\mathcal{D}
$$
where $s\in\mathcal{S}$ is the current state of the process and $D\in\mathcal{D}$ is the deterministic transition matrix that gives the state of the process on the next step. That way, a partial realization up to time $n$ and a complete realization of a RDS respectively take their values in the spaces
$$
\Omega_N =\hat{\Omega}^N = \left(\mathcal{S}\times\mathcal{D}\right)^N .
$$
\begin{equation}
\Omega =\hat{\Omega}^\mathbb{N} = \left(\mathcal{S}\times\mathcal{D}\right)^\mathbb{N} .
\label{eq:omegaRDS}
\end{equation}

%\textsl{`` Yet, a realization of a RDS is a sequence of matrices and a sequence of states, so we need to define
%$$
%\mathcal{D}^\infty = \left\{\left(D_n\right)_{n\in\mathbb{N}^*} : \forall n \in\mathbb{N}^*,  D_n \in\mathcal{D}\right\},
%$$
%$$
%\mathcal{S}^\infty = \left\{\left(s_i\right)_{i\in\mathbb{N}} : \forall i \in\mathbb{N},  s_i \in\mathcal{S}\right\}.
%$$
%Now we can define
%$$
%\Omega^\infty = \mathcal{D}^\infty \times \mathcal{S}^\infty
%$$ ``} Not sure about this definition.\\

Moreover, for the sake of simplicity, for each $\omega \in\Omega$, let us introduce $D_n(\omega) = D_n$ and $X_n(\omega) = s_n$, for any $n\in\mathbb{N}$. So, $X_0 = s_0$ is the initial state of the process, and $D_0$ is the deterministic transition matrix between step $0$ and step $1$.

From now on, we want to introduce a $\sigma$-field on $\hat{\Omega}$, $\Omega_n$ and $\Omega$ to properly define a probability space. The powersets $\mathcal{P(S)}$ and $\mathcal{P(D)}$ respectively are $\sigma$-fields for $\mathcal{D}$ and $\mathcal{S}$, so $\mathcal{P(S\times D)}$ is one for $\hat{\Omega}$. We can so define the $\sigma$-field $\mathcal{F}_N$ defined by the cylinder sets such that
$$
\mathcal{F}_N = \sigma \left(C^N_{n_1,\cdots,n_k,m_1,\cdots,m_l}(A_1,\cdots,A_k,B_1,\cdots,B_l)\right),
$$
with $(k,l)\in\left\{1,2,\cdots,N\right\}^2, A_i\in\mathcal{P(D)}, B_j\in\mathcal{P(S)}$
and
\begin{align*}
&C^N_{n_1,\cdots,n_k,m_1,\cdots,m_l}\left(A_1,\cdots,A_k,B_1,\cdots,B_l\right) =
\\
&\Big\{ \omega \in\Omega : D_{n_1}(\omega)\in A_1,\cdots,D_{n_k}(\omega)\in A_k
 X_{m_1}(\omega)\in B_1,\cdots,X_{m_l}(\omega)\in B_l \Big\}.
\end{align*}
$\left(\mathcal{F}_N\right)_{N\in\mathbb{N}^*}$ is a filtration on $\left(\Omega_N\right)_{N\in\mathbb{N}^*})$ and
$$
\mathcal{F} = \sigma \left(\bigcup_{N\in\mathbb{N}}C^N_{n_1,\cdots,n_k,m_1,\cdots,m_l}(A_1,\cdots,A_k,B_1,\cdots,B_l)\right)
$$
is a $\sigma$-field on $\Omega$.

Finally, let us consider the probability measure $\mathbb{P}^N_{\mu_0}$ (rep. $\mathbb{P}_{\mu_0}$) on $\left(\Omega_N,\mathcal{F}_N\right)$ (resp. $\left(\Omega,\mathcal{F}\right)$) defined by
\begin{align}
\forall \tilde{\omega}\in\Omega_N& : \nonumber\\
\mathbb{P}^N_{\mu_0}(\tilde{\omega}) &= \mu_0(s_0,D_0) \prod_{n\leq N}{q_n D_{n-1}(s_{n-1},s_n)} \label{eq:probaN1}\\
 &=\mu_0(s_0,D_0) \prod_{n\leq N}{q_n}\prod_{n\leq N}{ D_n(s_{n-1},s_n)} \label{eq:probaN2}\\
\forall \omega\in\Omega : &\nonumber\\
\mathbb{P}_{\mu_0}(\omega) &= \mu_0(s_0,D_0) \prod_{n\in\mathbb{N}^*}{q_n D_{n-1}(s_{n-1},s_n)} \label{eq:proba1}\\
 &=\mu_0(s_0,D_0) \prod_{n\in\mathbb{N}^*}{q_n}\prod_{n\in\mathbb{N}^*}{ D_n(s_{n-1},s_n)}
\label{eq:proba2}
\end{align}
where $\mu_0 : \hat{\Omega} \rightarrow \left[0,1\right]$ is the initial probability distribution on $\hat{\Omega}$. This is precisely the probability measure described later, which corresponds first to draw at each step a transition matrix, and then to effectuate the transition om the current state to the state indicated by the matrix. $D_n(s_{n-1},s_n)$ corresponds to the coefficient of the matrix $D_n$ associated to the transition from state $s_{n-1}$ to $s_n$. Here, this coefficient can only be either $0$ or $1$, due to the deterministic aspect of the problem. $0$ means that the transition $s_{n-1} \rightarrow s_n$ is not possible with the selected matrix, and so that the path $\left(s_0, s_1, s_2,\cdots\right)$ is impossible with this sequence of matrices.

Let us call $\Delta = \left(D_n\right)_{n\in\mathbb{N}}$, for $\omega\in\Omega$ the probability measure knowing $\Delta$ is
\begin{equation}
\mathbb{P}_{\mu_0}^{\Delta}(\omega) = \mathbb{P}_{\mu_0}(\omega|\Delta) = \mu_0(s_0|D_0) \prod_{n\in\mathbb{N}^*}{D_n(s_{n-1},s_n)}.
\label{eq:probadelta}
\end{equation}
We notice that, for a RDS, when the starting state is known, $\mathbb{P}_{\mu_0}^{\Delta}(\omega) \in \left\{0,1\right\}$, which means that, once the transition matrices are drawn, there exists only one path possible between the state. This is the deterministic idea. The knowledge of the sequence of transition matrices brings the knowledge of the sequence of states (however, the opposite is rarely true).

\begin{figure}[!h]%
\includegraphics[width=\columnwidth]{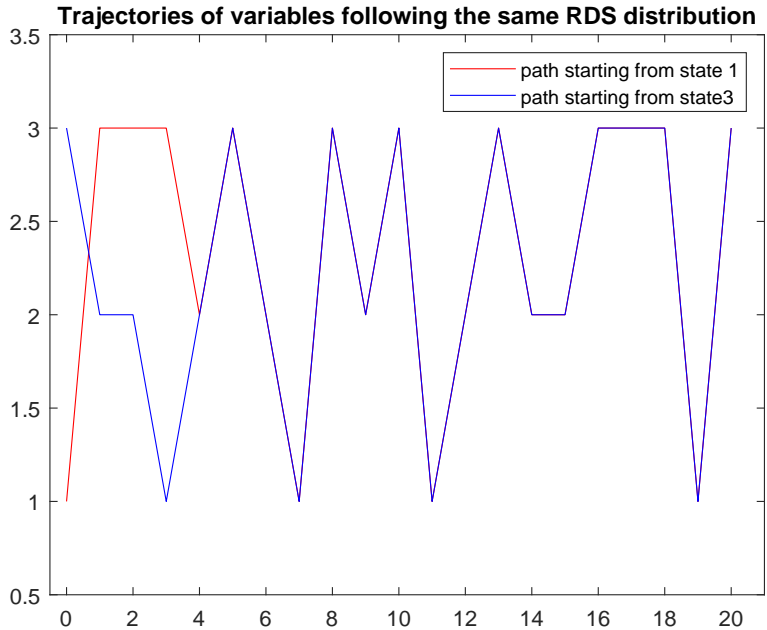}%
\caption{20 first steps of a paths following the same RDS distribution, starting from different states. We see here a synchronization.}%
\label{fig:rds3}%
\end{figure}

\begin{figure}[!h]%
\includegraphics[width=\columnwidth]{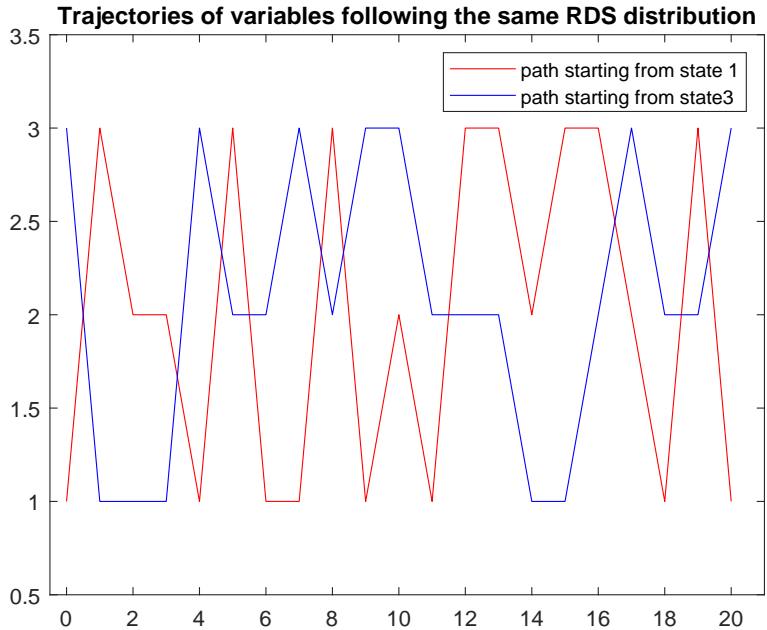}%
\caption{With this set of matrices, there is no synchronization possible.}%
\label{rds31}%
\end{figure}

It is interesting to also notice that this description of the RDS allow us to study the synchronization issue of two paths following the same RDS i.i.d. distribution. Indeed, for $\omega_1,\omega_2 \in\Omega$, conditioning on $\Delta(\omega_1) = \Delta(\omega_2)$ means that the two sequences of states can be seen as drawn simultaneously with the same transition matrices. Originally, the purpose of the RDSs is to induce the same transition's behavior to all the entities of a single system, so it is important to keep this application in mind. That is exactly why one uses RDSs as a model of extrinsic noise. We can see the phenomenon of synchronization in Figure \ref{fig:rds3}. However, some sets of parameters $\left(q_i,D_i\right)$ do not lead to synchronization (see Figure \ref{rds31}).

\subsection{Perturbation into a Random Markov System (RMS)}
\label{RMS}
The goal now is to introduce a second type of noise during the transitions from one state to another, while keeping the extrinsic noise. This noise susceptible to modify a single transition without acting on the other entities' one is called intrinsic noise. So we want to keep the structure of the RDS, adding in the fact that at each step there exists a very small probability not to follow the transition of the deterministic matrix, but another one.

Let $Q$ be a matrix of size $s\times s$ such that each row sums to $0$, with $-1$ on the diagonal coefficients, and $\epsilon \in \mathbb{R}_+$ a non-negative real number. We notice that the matrices $Q^n$ for any $n\in\mathbb{N}^*$ respect also the same property, the sum over each row is still equal to $0$, and so $N = e^{\epsilon Q}=exp(\epsilon Q)$ is a Makovian matrix. If we consider a deterministic transition matrix $D\in\mathcal{D}$, the product $D N = D e^{\epsilon Q}$ gives the matrix with repositioned rows. We will discuss this more precisely later, but we can already notice that, when $\epsilon$ gets smaller, $e^{\epsilon Q} \approx_{\epsilon \rightarrow 0} \left(I_s-\epsilon Q\right)$, which diagonal coefficients are $1-\epsilon$. So two noteworthy facts are happening. First, in each row, the highest coefficient of $D\times N$ has the same position that the $1$ in $D$, so $DN$ is still close to $D$. Second, by comparing $D$ and $D N$, $\epsilon$ can be interpreted as the probability of not following the deterministic transition of $D$, as long as $\epsilon$ is small enough.

As in \ref{RDS}, we can easily introduce a well defined probability space using the same $\Omega$ and $\sigma$-field $\mathcal{F}$, and substituting $D_n(s_{n-1},s_n)$ by $(D_n N)(s_{n-1},s_n)$ in \eqref{eq:proba1}. This description is equivalent as describing a Markovian process in which we draw at each step a Markov transition matrix. These kind of processes are called Markov Processes with Random Transition Matrix, or more easily Random Markov Systems (RMS). In our situation, we draw the matrix $D_i N$ with the probability $\left(q_i\right)$, where $D_i\in \mathcal{D}$. This way we turn the previous RDS into a RMS, with an adequate mathematical environment to study it. However, the interesting point is to study at the same time both the RDS and the RMS. This would allow us to measure to what extent they are dependent to each other, and to what extent this dependency vary with $\epsilon$. We shall so introduce a new probability space that record at each step the deterministic transition matrix $D$ drawn, the state of a path following the deterministic transition (RDS) and the state of another path following transition which probability are given by $D N$ (RMS). We will note $X = \left(X_n\right)\in\mathcal{S}^\mathbb{N}$ the path from the RDS, and $Y = \left(Y_n\right)\in\mathcal{S}^\mathbb{N}$ the one from the RMS.

We define
$$
\hat{\Omega} = \left(\mathcal{S} \times \mathcal{S} \times \mathcal{D}\right)
$$
$$
\Omega_N = \hat{\Omega}^N = \left(\mathcal{S} \times \mathcal{S} \times \mathcal{D}\right)^N
$$
\begin{equation}
\Omega = \hat{\Omega}^\mathbb{N} = \left(\mathcal{S} \times \mathcal{S} \times \mathcal{D}\right)^\mathbb{N}
\label{eq:omegaRMS}
\end{equation}
and, as we did before, we can define the new $\sigma$-fields $\mathcal{F}_N$ and $\mathcal{F}$ using the cylinder sets (with three sequences of indexes, instead of two). The interesting part is to define the probability of a joint transition of $X$ and $Y$. With $\mu_0$ as the initial distribution probability on $\hat{\Omega}$. For $\omega\in\Omega$, we have:
%\begin{align}
%\mathbb{P}_{\mu_0}(\omega) &= \mu_0(x_0,y_0) \prod_{n\in\mathbb{N}^*}{q_n D_n(x_{n-1},x_n) D_n N(y_{n-1},y_n)} \label{eq:proba_joint1} \\
 %&=\mu_0(x_0,y_0) \prod_{n\in\mathbb{N}^*}{q_n}\prod_{n\in\mathbb{N}^*}{ D_n(x_{n-1},x_n) D_n N(y_{n-1},y_n)}
%\label{eq:proba_joint2}
%\end{align}

\begin{align}
\forall \tilde{\omega}\in\Omega_N& : \nonumber \\
\mathbb{P}^N_{\mu_0}(\tilde{\omega}) &= \mu_0(x_0,y_0,D_0) \prod_{n\leq N}{q_n D_{n-1}(x_{n-1},x_n) D_{n-1}N(y_{n-1},y_n)} \label{eq:proba_jointN1}\\
 &=\mu_0(x_0,y_0,D_0) \prod_{n\leq N}{q_n}\prod_{n\leq N}{ D_{n-1}(x_{n-1},x_n) D_{n-1}N(y_{n-1},y_n} \label{eq:proba_jointN2}\\
\forall \omega\in\Omega : &\nonumber\\
\mathbb{P}_{\mu_0}(\omega) &= \mu_0(x_0,y_0,D_0) \prod_{n\in\mathbb{N}^*}{q_n D_{n-1}(x_{n-1},x_n) D_{n-1}N(y_{n-1},y_n)} \label{eq:proba_joint1}\\
 &=\mu_0(x_0,y_0,D_0) \prod_{n\in\mathbb{N}^*}{q_n}\prod_{n\leq N}{ D_n(x_{n-1},x_n)D_{n-1}N(y_{n-1},y_n}
\label{eq:proba_joint2}
\end{align}

Note that we carefully draw the same matrix $D_n$ at each step for both the paths $X$ and $Y$, but they don't have the same transition probability thanks to $N$.

\begin{figure}[!h]%
\includegraphics[width=\columnwidth]{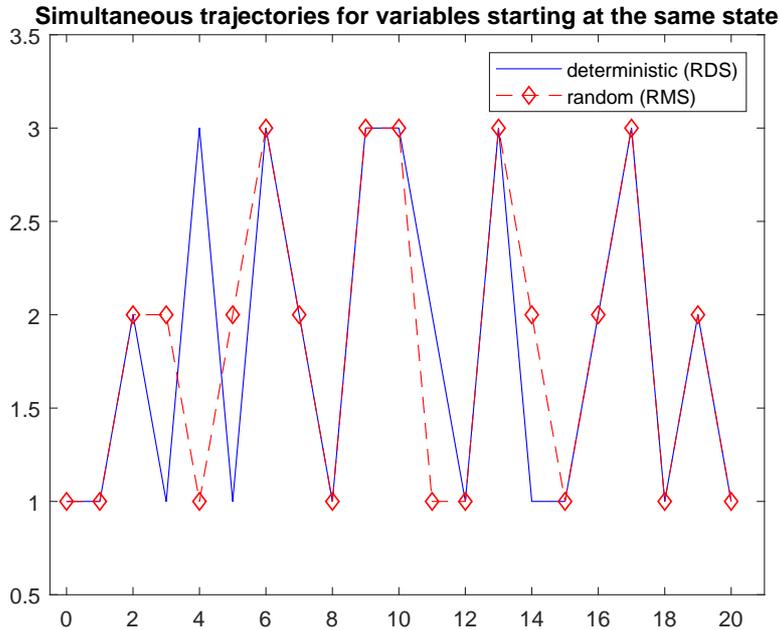}%
\caption{20 first steps of the two paths described above, with $\epsilon = 0.1$ and a random $Q$.}%
\label{fig:rms3}%
\end{figure}

In the end, we only need to know the sequences of states taken by the RDS and the RMS, however we need the sequence of matrices to compute the probabilities. That's why we need to keep $\left(D_n\right)$ in each $\omega$. When the RDS and the RMS are in the same state, the probability $\epsilon$ that they are not together in the same state is very small, so we should see them go together for quite a long time. Then, since this probability is non zero, at some point they will diverge, or desynchronize. From there, they are nearly behaving as two variables following the same RDS, since the perturbation matrix $N$ is close from $I_s$, and so we might see a phenomenon of synchronization happening between $X$ and $Y$. Finally, at this point, we can think that the process got refresh, we could expect them to behave again as if they were starting from the same initial state. We can notice this in the Figure \ref{fig:rms3}.

\subsection{Properties and approximations}
\label{approx}

Mathematically, we have the very simple but essential result, as defined in section \ref{RMS}, the process of comparing the RDS and the RMS is a Markov Chain (MC).
\begin{prop}
The process $\left(H_n\right)_{n\in\mathbb{N}}$ defined by
$$
H_n(\omega) = \left( X_n(\omega), Y_n(\omega), D_n(\omega)\right)
$$
for $n\in\mathbb{N}$, $\omega\in\Omega$ is a Markov Chain regarding the filtration $\left(\mathcal{F}_n\right)$.
\end{prop}
This relies on the equality
\begin{align} 
\mathbb{P}_{\mu_0}\left(H_n=h_n\middle| H_0=h_0, H_1=h_1,\cdots,H_{n-1}=h_{n-1} \right) &= q_n D_{n-1}(x_{n-1},x_n) D_{n-1}N(y_{n-1},y_n) \nonumber \\
&= \mathbb{P}_{\mu_0}\left(H_n=h_n\middle| H_{n-1}=h_{n-1}\right)
\label{eq:rmsmc} 
\end{align}
for any $h_i=\left(x_i,y_i,D_i\right)\in\hat{\Omega}$.\\ Obviously, the process defined by the RDS alone is also a MC.

Let us discuss a bit longer about the part of $\epsilon$ and $Q$ in the calculus of $N$. We want to understand what value we have to take for these two parameters, in order to model a biological behavior, or at least a plausible noise.\\

\begin{figure}%
\includegraphics[width=\columnwidth]{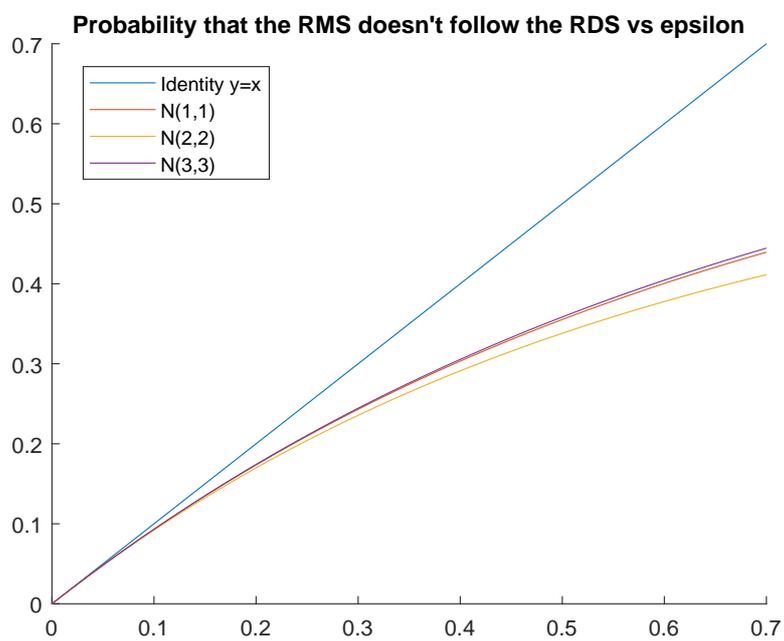}%
\caption{We have plot the diagonal coefficients of $N$ according to $\epsilon$, with a random $Q$. We notice that the gap between the blue line and the others is negligible under $0.1$, but starts getting wider after. Here, the size of $Q$ $s=3$. We have similar results with $s$ higher, the line are less scattered, which means that the diagonal coefficients are closer from each other.}%
\label{variationQeps}%
\end{figure}

Let first focus on $\epsilon$. As we already said in the previous section \ref{RMS}, as $\epsilon$ gets smaller, $e^{\epsilon Q}$ gets closer to $I_s-\epsilon Q$. The Taylor expansion of the matrix exponential shows us that the gap between $e^{\epsilon Q}$ and this very approximation is in the range of $\epsilon^2$. So, when $\epsilon^2$ can be neglected regarding to $\epsilon$, then $I_s-\epsilon Q$ becomes a very good approximation of this noise. As shown in figure \ref{variationQeps}, we shall keep in mind that $0<\epsilon<0.1$ is a good area where this approximation holds. In practice, we remember that small $\epsilon$ is the probability that, at some point, the trajectory of the RMS doesn't follow the one given by the RDS. $\epsilon$ is so a sort of amplitude of the noise $N$ that impacts the RDS, and therefore it's reasonable to imagine that it's (way) smaller than $0.1$ for practical applications. We this probability-oriented understanding, we notice that the probability $\epsilon$ not to follow the RDS is split between the other transitions possible, and so that the non-diagonal coefficients of each row of $N$ should sum to something close to $\epsilon$.\\

\begin{figure}[h!]%
\includegraphics[width=\columnwidth]{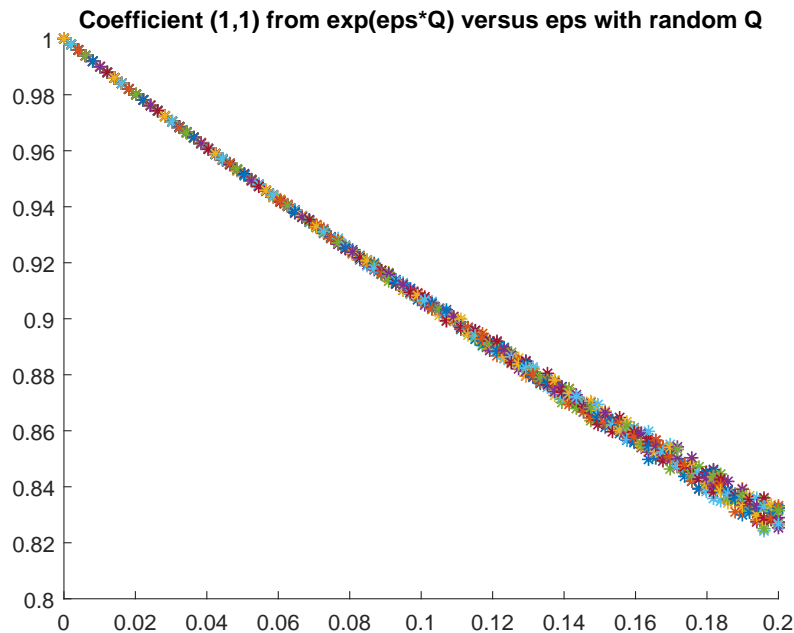}%
\caption{For each value of $\epsilon$, we have plot 10 realizations of Q. he number of states is s=3 (dimension of the square matrix $Q$).}%
\label{variabilityQ}%
\end{figure}

Let us now take a closer look to the role of $Q$, especially its diagonal coefficients. When $s=2$, $Q$ is completely determined by its definition, and we have
\begin{equation}
Q = \left( \begin{array}{cc}
-1 & 1\\
	1 & -1\\
\end{array}\right)
\label{eq:Q2}
\end{equation}

and so it is very easy to compute (by diagonalizing)
$$
N=e^{\epsilon Q} = \frac{1}{2}\left( \begin{array}{cc}
1+e^{-2\epsilon} & 1-e^{-2\epsilon}\\
	1-e^{-2\epsilon} & 1+e^{-2\epsilon}\\
\end{array}\right).
$$
When $s=3$, $Q$ has the shape
$$
Q = \left( \begin{array}{ccc}
-1 & a & b\\
c & -1 & d\\
e & f & -1\\
\end{array}\right)
$$
with
$$
\left\{ \begin{array}{c}
	a+b=1\\
	c+d=1\\
	e+f=1\\
\end{array} \right. .
$$
Analytically, it is still possible to calculate $N$ though it is not very relevant here. Indeed, even if we have emphasized some quantities that give information on how $N$ varies with $Q$ (and its coefficients), this variation is not significant. The figure \ref{variabilityQ} shows how the first coefficient from $N$ evolves. For each value of $\epsilon$, we draw $10$ different values of $Q$, each non-diagonal coefficient following an uniform distribution over $\left[0,1\right]$. As $\epsilon$ increases, we notice that the data are more scattered, but the vertical expansion of the values is not significant compared to the change induced by $\epsilon$ (derivative with respect to $\epsilon$). This gets even more true as $\epsilon$ gets smaller, so well that we can neglect this $Q$-dependance under $0.1$. \\

\begin{figure}[h!]%
\includegraphics[width=\columnwidth]{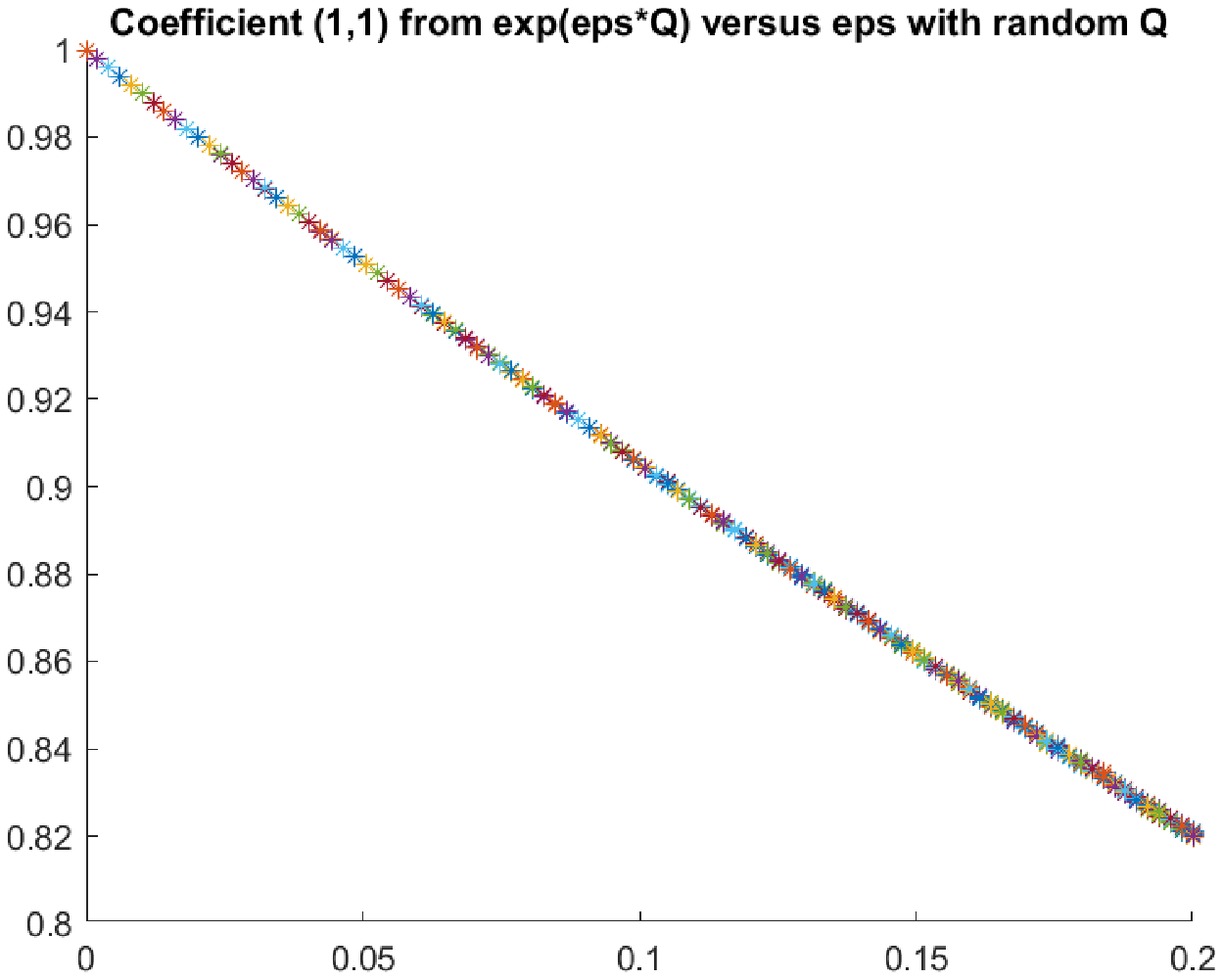}%
\caption{For each $\epsilon$, 100 values of $Q$ are drawn, where $Q$ has the size $s=12$.}%
\label{variabilityQ12}%
\end{figure}

When $s$ is higher, it is not humanely reasonable to calculate the expression of $N$ according of the coefficients of $Q$. Even if it were, it would not be helpful insofar as it would probably depend on a lot of relations between the coefficients of $Q$, which number increases really fast. Hopefully, numerical computation shows us that the increase of $s$ reduces such vertical expansion, inasmuch as the randomly distributed $Q$ matrices transform into $N=e^{\epsilon Q}$ matrices that are more and more condensed (as we can see in Figure \ref{variabilityQ12}).

However, even if the non-diagonal coefficients of $Q$ do not really influence the diagonal values of $N$, they are very sensitive on the rest of the matrix. So, keeping a random $Q$ (all of its non-diagonal coefficients are drawn with a uniform distribution on $\left[0,1\right]$, then normalized row by row), contributes to keep all of the non-diagonal coefficients of $N$ very random from one realization to another. Nevertheless, we can wonder about the relevance of keeping the same Q during all the process. For a deterministic transition matrix $D$ from the RDS, if we go into details about the shape of $D N$, for each starting starting state of the transition, there is a high probability $1-\epsilon$ to do the transition given by $D$ for this starting state, and then a probability $\epsilon$, split into all of the other states, to do the transition from the same starting state to those others. Now, if we keep the same starting state, but we change $D$. Then the arriving state given by $D$ might change, which means that another row from $N$ could be used as the transition probability from the same starting state (we remind that the product $D\times N$ simply exchange the rows of $N$ amongst themselves), and so the noise  distribution from this very starting state would be different. We realize that the noise (probability) distribution of the non-favored transitions does not depend on the starting state of the transition, but on the arriving state. Consequently, with a non-time-dependent $Q$ (and so on for $N$), we might notice on a long time scale some repetitive transition behavior of bifurcation from the RDS. Indeed, if, according to the RDS, a transition was supposed to happen towards a given state, the probability distribution of the transitions towards the other states are the same along the time, independently of the starting state of this transition. And so, knowing that a bifurcation between the RDS and the RMS is supposed to happen, independently of the time at which this transition happen, we know the more likely states for this bifurcation to arrive to. This wouldn't be true with a changing $Q$ along the time, but, as we have shown, the probability no to follow the RDS would still be similar (for $\epsilon$ small enough).

\section{Relative study of the RMS: Synchronization rate}
\label{synchro}

Now, we are starting to have a good understanding of how this process works and evolves along the time. We have seen what is a RDS, how we can turn it into a RMS, and to what extent this transformation happens. We also have the appropriate mathematical environment to process to an analytical study.

The leading phenomenon that seems to be at the root of the long time behavior is a repetition pattern pattern. Indeed, let's start with the RDS and the RMS at the same state. They will evolve together until a desynchronization happens (due to the RMS), and then they will go each one on their own way to the next synchronization. At this point, we feel like the process has refreshed, and is starting again (possibly from another state than previously, but we actually focus on the two states synchronized/unsynchronized problem). Some questions naturally emerge from this reasoning. What is the rate of synchronized time over the total time of the process? Do we have the existence of a mean or pseudo period that would characterize the process? The renewal process theory can help us address these questions.

\subsection{Random Times}
\label{times}
We now want to clearly decompose the sequences of states pointing out the interesting times when the RDS and the RMS synchronize or desynchronize. Let introduce some stopping times on $\Omega$. The times $T_i$ will count the duration of the successive synchronization phases, from the $i-1$-th first synchronized state to the next $i$-th first one after a desynchronization. $\tau_i$ will measure the duration of the $i$-th unsynchronized phase. We also introduce the waiting time until the $i$-th synchronization $W_i$ as the sum over all the previous $T_k$ for $k\leq i$. Mathematically,  for any $\omega\in\Omega$, we have

\begin{align}
&T_0(\omega) = 0 \label{T0},\\
&\forall i\in\mathbb{N}^*, \left\{\begin{array}{c}
W_i = \sum_{k=0}^i{T_k} \\
\tau_i(\omega) = min\left\{n\in\mathbb{N}^* : X_{W_{i-1}+n}(\omega) \neq Y_{W_{i-1}+n}(\omega)\right\} \\
T_i(\omega) = min\left\{n>\tau_i : X_{W_{i-1}+n}(\omega) = Y_{W_{i-1}+n}(\omega)\right\} \\
\end{array} \right.
\label{stoppingtimes}
\end{align}
Since we will need it later, let also introduce the strict synchronization time $\gamma_i$ as the time needed from the $i$-th first desynchronized state to the $i$-th synchronized one defined by
\begin{equation}
\gamma_i = T_i - \tau_i
\label{eq:gamma}
\end{equation}
For example, if we report on the figure \ref{fig:rms3}, we have
$$
T_0 = 0 ,\  \tau_1 = 3 ,\  T_1 = 6 ,\  \tau_2 = 5 ,\  T_2 = 6 ,\  \tau_3 = 2 ,\  T_3 = 3 ,\  \tau_4 > 5 
$$
and, more generally
$$
\begin{array}{cc}
	X_{\tau_1} \neq Y_{\tau_1} & X_{T_1} = Y_{T_1}\\
	X_{W_1+\tau_2} \neq Y_{W_1+\tau_2} & X_{W_2} = Y_{W_2}\\
	X_{W_2+\tau_3} \neq Y_{W_2+\tau_3} & X_{W_3} = Y_{W_3}\\	
\end{array}
$$

Given all of this, we can imagine that an interesting number to quantify the process would be the rate of synchronized time for the RDS and the RMS over the total time of the experiment. It seems natural that the times $T_i$ are mutually independent, as well as the couples $\left(T_i,\tau_i\right)$. So on, one can imagine that the process is a repetition of independent and identically distributed cycle, so that the interesting ratio would be close to the expected values of $\tau$ over $T$. The renewal process theory specifically focuses on such problems, and shall bring answers and proofs.

\subsection{Double path transition matrix}
\label{transition_mat}

As we have already said, the mean transition matrix for a single path process $M = \sum_{i\in\mathcal{I}}{q_iD_i}$ only catches the behavior of the MC corresponding to the RDS (same for the RMS with $\hat{M} = \sum_{i\in\mathcal{I}}{q_iD_iN}=MN$). However, if we want to study synchronizations (and so on the times $T_i$ and $\tau_i$), we need to focus on transitions on a ``double paths'' process, which would consider the states of two simultaneous process. Before studying the RDS/RMS comparison, let us first explain this idea on two paths following the same RDS distribution, as given in section \ref{RDS}.

As before, using the same space defined in \eqref{eq:omegaRDS}, we get into the condition $\Delta(\omega_1)=\Delta(\omega_2)$, for any $\omega_1,\omega_2\in\Omega$. Let call $(u_i)$ and $(v_i)$ the sequences of states of the realizations $\omega_1$ and $\omega_2$. The states of the double path process are couples of states from $\mathcal{S}$. If, for some time $i\geq 0$, $u_i=v_i$, then the paths will stay together forever, for any $j\geq i : u_j=v_j$, and so the mean transition probability $\mathbb{P}\left(u_{j+1}\middle|u_j\right) = M(u_{j+1},u_j)=\sum_{i\in\mathcal{I}}{q_iD_i(u_{j+1},u_j)}$. It's actually the mean transition probability for a single path RDS process. But when the states $u_i$ and $v_i$ are different, we need another tool to describe the possible transitions. Let call $V$ the mean transition matrix of the double path process (from $\mathcal{S}^2$ to itself) we have 
$$V_{(u_j,v_j)\rightarrow(u_{j+1},v_{j+1})}=\mathbb{P}\left(u_{j+1},v_{j+1}\middle|u_j,v_j\right) = \sum_{i\in\mathcal{I}}{q_iD_i(u_j,u_{j+1})D_i(v_j,v_{j+1})}.
$$
 The size of $V$ is so $s^2\times s^2$. If we rearrange the rows and the line of $V$ in order to have the states $(1,1),(2,2),\cdots,(s,s)$ in the last $s$ rows and $s$ columns, then $V$ has an interesting shape. We have to notice that $V$ is still a  Markovian matrix.
\begin{equation}
V=\left(
\begin{array}{c|c}
\begin{array}{c c c c}
	\dots & \dots & \dots & \dots \\
	\vdots & & & \vdots \\
	 & & & \\
	\vdots & & & \vdots \\
	\dots & & & \dots \\	
\end{array}
&
\begin{array}{c c}
	\dots & \dots \\
	\vdots & \vdots \\
	 & \\
	\vdots & \vdots \\
	\dots & \dots \\
\end{array} \\
\hline
 & \\
(0) & M	
\end{array}
\right)
=
\left(
\begin{array}{c|c}
\bar{M} &
\begin{array}{c c}
	\dots & \dots \\
	\vdots & \vdots \\
	 & \\
	\vdots & \vdots \\
	\dots & \dots \\
\end{array} \\
\hline
 & \\
(0) & M	
\end{array}
\right) 
\label{eq:matV}
\end{equation}

We can distinguish the states of $\mathcal{S}^2$ between the synchronized states (of the shape $(i,i)$) and the others ( $(i,j)$ with $i\neq j$). In the upper left square, in the matrix $\bar{M}$, there are the mean probabilities to go from an unsynchronized state to another one, and the lower right square shows the mean probabilities to stay in a synchronized state coming from an already synchronized one. More interestingly, the upper right rectangle displays the mean probabilities to synchronize, and the lower left one to unsynchronize. Obviously, for two RDS distributed paths, the latter is equal to zero : it's impossible given the same deterministic transition matrix to go from two instances of a same state to two different states. Let us use the matrix given in \eqref{eq:mean_markov} as an example. The first decomposition of this matrix gives the double paths mean transition matrix
$$
V_1 = \left( \begin{array}{cc|cc}
	0.2 & 0.6 & 0 & 0.2\\
	0.6 & 0.2 & 0 & 0.2\\
	\hline
	0 & 0 & 0.2 & 0.8\\
	0 & 0 & 0.6 & 0.4\\
\end{array}\right)
$$
and the second one gives
$$
V_2 = \left( \begin{array}{cc|cc}
	0.1 & 0.5 & 0.1 & 0.3\\
	0.5 & 0.1 & 0.1 & 0.3\\
	\hline
	0 & 0 & 0.2 & 0.8\\
	0 & 0 & 0.6 & 0.4\\
\end{array}\right).
$$

Given the matrix $V$, we can easily calculate the expected value of the first synchronization time from an unsynchronized state, let us call it $\gamma_1$ (the notation is not rigorous, as $\gamma_i$ has been defined for the double paths process RDS/RMS only). Let $\mu_0$ be the initial distribution over $\mathcal{S}^2$, 
\begin{equation}
\begin{split}
\mathbb{E}_{\mu_0}\left(\gamma_1\right)&=\sum_{n=1}^{\infty}{\mathbb{P}_{\mu_0}\left(\gamma_1\geq n\right)}\\
&=\sum_{n=1}^{+\infty}{\sum_{\begin{subarray}{c}
	(i,j)\in\mathcal{S}^2\\
	i \neq j
\end{subarray}}{\mu_0(i,j)\mathbb{P}_{\left(i,j\right)}{\left(\gamma_1\geq n\right)}}} \\
&=\sum_{\begin{subarray}{c}
	(i,j)\in\mathcal{S}^2\\
	i \neq j
\end{subarray}} {\mu_0(i,j)\sum_{n=1}^{+\infty}{\mathbb{P}_{\left(i,j\right)}\left(u_0\neq v_0,u_1\neq v_1,\cdots,u_{n-1}\neq v_{n-1}\right)}}
\end{split}
\label{eq:expected_value_T1}
\end{equation}
Yet, $\mathbb{P}_{\left(i,j\right)}\left(u_0\neq v_0,u_1\neq v_1,\cdots,u_{n-1}\neq v_{n-1}\right)$ is pretty easy to calculate. It is the probability to go in $n-1$ steps from the state $(i,j)$ to any other state $(k,l)$ with $k\neq l$, without going through any state $(m,m)$ in between. If we note $\mathbb{I}$ the column vector of size $s(s-1)$ with only ones, then we have
\begin{equation}
\mathbb{P}_{\left(i,j\right)}\left(\gamma_1\geq n\right) = \bar{M}^{n-1}\cdot\mathbb{I},
\label{eq:pij}
\end{equation}
and so, if we note $\mu$ the row vector of the values $\mu_0(i,j)$ for $i\neq j$, then we obtain the nice equation

\begin{equation}
\mathbb{E}_{\mu_0}\left(\gamma_1\right) = \sum_{n=1}^{+\infty}{\mu_0\cdot \bar{M}^{n-1}\cdot\mathbb{I}} = \mu_0\cdot\left( \sum_{n=1}^{+\infty}{\bar{M}^{n-1}}\right)\cdot\mathbb{I}.
\label{eq:expT1}
\end{equation}

Finally, we can apply a last modification to $V$. One may notice that we don't really pay attention the the precise state of the system but rather on the fact that it is synchronized or not. So, it makes sense to gather all the synchronized states in only one, let us call it $S = \left\{(i,i): i\in\mathcal{S}\right\}$, and now gather the last $s$ columns of $V$ into only $1$ by adding them together. Each coefficient of the last column shall be the probability to go from  the state represented by the row to one of any state from $S$. Then, we also gather the last $s$ rows into a single $1$, according to the fact that the probability to go from a synchronized state to an unsynchronized one is $0$, and $1$ to go from $S$ to $S$. Let call this new matrix $\tilde{V}$. With the previous examples, we obtain
$$
\tilde{V}_1 = \left( \begin{array}{cc|c}
	0.2 & 0.6 & 0.2\\
	0.6 & 0.2 & 0.2\\
	\hline
	0 & 0 & 1\\
\end{array}\right)
$$
and
$$
\tilde{V}_2 = \left( \begin{array}{cc|c}
	0.1 & 0.5 & 0.4\\
	0.5 & 0.1 & 0.4\\
	\hline
	0 & 0 & 1\\
\end{array}\right).
$$
\newline

Now that we have seen the basics, let us get into what really matters for us. We want to introduce the matrices $W$ and $\tilde{W}$ in the same way, not to compare two RDS simultaneous paths, but a RDS and a RMS ones, as presented in our model from section \ref{RMS}. With the same notation, we can write 
$$
W_{(x_j,y_j)\rightarrow(x_{j+1},y_{j+1})}=\mathbb{P}\left(x_{j+1},y_{j+1}\middle|x_j,y_j\right) = \sum_{i\in\mathcal{I}}{q_iD_i(x_j,x_{j+1})D_iN(y_j,y_{j+1})}.
$$
In order to make it more readable, we shall rearrange $W$ the same way we did with $V$, keeping the states from $S$ on the last $s$ rows and columns. $W$ is now like $V$ in \eqref{eq:matV}, except that the lower left rectangle is no longer empty. Indeed, the probability to desynchronize (to go out $S$) is not equal to $0$ but depends on $N$. More explicitly, it depends on the diagonal coefficients of $N$. If we still call $M$ the mean value of the transition matrix $M=\sum_{i\in\mathcal{I}}{D_iN}$, then the lower right square of $W$ is equal to $M\cdot Diag\left(N\right)$ where $Diag\left(N\right)$ is the diagonal matrix of the diagonal coefficients of $N$. As we have already studied in section \ref{approx}, when $\epsilon$ gets smaller, $Diag\left(N\right)$ converges towards $(1-\epsilon)I_s$, and so the sum of the coefficients of each row of the lower right square of $W$ converges towards $1-\epsilon$. We finally have the expected result that the unsynchronization probability is close to $\epsilon$, which can be found by summing over the columns on each row of the lower left rectangle of $W$.

We notice that, given $\epsilon$, $Q$ and the probabilities $(q_i)$, we can calculate in the same time $V$ and $W$. Let us so do it for $s=2$ with a random probability distribution $(q_i)$ for $i\in\left\{1,2,3,4\right\}$ and $\epsilon=0.01$. A short MATLAB program gives
$$
V = \left( \begin{array}{cc|cc}
	0.1588 & 0.4759 & 0.3480 & 0.0173\\
	0.4759 & 0.1588 & 0.3480 & 0.0173\\
	\hline
	0 & 0 & 0.5068 & 0.4932\\
	0 & 0 & 0.8239 & 0.1761\\
\end{array}\right)
$$
and
$$
W = \left( \begin{array}{cc|cc}
	0.1607 & 0.4714 & 0.3461 & 0.0218\\
	0.4747 & 0.1574 & 0.3493 & 0.0187\\
	\hline
	0.0050 & 0.0049 & 0.5018 & 0.4883\\
	0.0082 & 0.0017 & 0.8158 & 0.1743\\
\end{array}\right).
$$
The formula \eqref{eq:expT1} still holds for $W$, using $\bar{M}$ as the upper left square of size $s(s-1)$ and $\mu_0$ as the distribution probability over the unsynchronized states, knowing that the initial state is not in $S$. However, a new issue arises here. With the RDS/RDS comparison process, at most one synchronization may happen, and so we only need once the initial distribution, but for the RDS/RMS process, there can be an infinite number of desynchronization, and so we wonder which probability distribution should be used in the formula as the initial distribution of each cycle (the $i$-th cycle is everything that happens between the times $W_{i-1}$ and $W_i$, of duration $T_i$). In order to tackle this question, we need to focus on what happens during the $\tau_i$ phase, when the RDS and the RMS paths are synchronized.

Let imagine that, at time $i$, $X_i=Y_i$, then the probability that they are still equal at time $i+1$ would be the coefficient at the position $(X_{i+1},X_{i+1})$ in N. Yet, when $\epsilon$ gets small enough, all the coefficients $N(k,k)$ get closer from $1-\epsilon$, and so this probability tends not to depend any longer from $X_{i+1}$ as all the coefficients have the same limit. Then the probability of desynchronization after $1$ transition is very close to $\epsilon$, after $2$ transitions is $(1-\epsilon)\epsilon$ and so on, such that $\mathbb{P}\left(X_{i+1}=Y_{i+1},\cdots,X_{i+n}=Y_{i+},X_{i+n+1}\neq Y_{i+n+1}\middle|X_{i}=Y_{i}\right)\approx \left(1-\epsilon\right)^n\epsilon$ . So, each $\tau_i$ is  very close to follow a geometric distribution on $\mathbb{N}^*$ with success probability $\epsilon$. That way, we can approximate a mean expected value for any $\tau_i$ with $\mathbb{E}(\tau_i) \approx \frac{1}{\epsilon}$. So, on average, the process will do $\frac{1}{\epsilon}-1 \approx \frac{1}{\epsilon}$ transitions according the RDS deterministic transition matrices before desynchronizing.

During this time, both motions behave as the same single path MC with the mean matrix $M$ as transition matrix, and so we can expect that the probability distribution between the states of $S$ will converge towards the steady state distribution given by $M$. As this convergence is exponentially fast, according Perron-Frobenius theorem (we can easily assume that the condition of the theorem will usually be satisfied, for example if we use a random probability set $(q_i)$). Given this reasoning, it sounds fair to estimate that the synchronized paths have reached the steady state distribution of the mean MC at time of desynchronization $\tau_i$. This remark is essential as we want to transform $W$ into $\tilde{W}$ as we did with $V$. Let us call $\pi(i)$ the steady state distribution of the Markovian matrix M, then the mean probability to go from $S$ to a state $(i,j)$ with $i\neq j$ will be 
\begin{equation}
\mathbb{P}\left(S\rightarrow(i,j)\right) = \sum_{l\in\mathcal{S}}{\pi(l)W_{(l,l)\rightarrow(i,j)}}.
\label{eq:proba_desynch}
\end{equation}
Using this formula, we can merge the last $s$ rows of $W$ into $1$, each coefficient being the mean unsynchronization probability that lead to the states $(i,j)$ (given by \eqref{eq:proba_desynch}). The last coefficient will be the mean probability to stay in $S$, which is supposed to be very close from $1-\epsilon$. With the previous example, we obtain
$$
\tilde{W} = \left( \begin{array}{cc|c}
	0.1607 & 0.4714 & 0.3679\\
	0.4747 & 0.1574 & 0.3679\\
	\hline
	0.0062 & 0.0037 & 0.9901\\
\end{array}\right).
$$
Let us call $\mu_1$ the probability distribution over the unsynchronized states knowing that there had been a desynchronization. $\mu_1$ is so the normalized last row of $\tilde{W}$ without the last coefficient (its size is $s(s-1)$). $\mu_1$ is also the mean ``initial distribution'' of any re-synchronization phase, it means that it is the probability measure over the states of $\left(\mathcal{S}^2\setminus S\right)$ at any time $W_i+\tau_i$. Finally, we can still use the formula \eqref{eq:expT1} to calculate the mean value of $\gamma_i$.\\
\newline
All the results discussed above can be easily guessed by a (somewhat) careful reading of the process. Let us try know to explain these results with more regressed maths. Let introduce the shifting (delay) operator $\theta$ on the process $(H_n)$ such that
\begin{equation}
\left\{\begin{array}{c}
	H_i\circ\theta_1 = H_i\circ\theta = H_{i+1}\\
	H_i\circ\theta_n = H_i\circ\theta\circ\theta_{n-1} = H_{i+1}\circ\theta_{n-1}
\end{array}\right. .
\label{eq:theta}
\end{equation} 
We want to prove that the random times $(\tau_i)$ and $(\gamma_i)$ are respectively independent and identically distributed (and so on same for $T_i$), on average. By ``on average'', we mean by marginalizing over the sequence of transition matrices $(D_i)$. As we did with the calculus of $W$, we want to know the mean behavior of the process, so we can consider the simultaneous double point motion $(X_n,Y_n)$ as a MC with transition matrix $W$. With the approximation of small $\epsilon$, each $\tau_i$ follows a geometric law, and so, even if they are unbounded, they are finite almost surely. With our model of the intrinsic noise $N$ as an exponential matrix, the rightest column of $\tilde{W}$ is non equal to only zeros, and so the probability that any $\gamma_i=\infty$ is $0$. So, for any $i\in\mathbb{N}$, we have
\begin{equation}
\mathbb{I}_{\left\{W_i<\infty\right\}} = 1 \text{     a.s.}
\label{eq:wi}
\end{equation}
One can also notice the important following relations
\begin{align}
\tau_i = \tau_1\circ \theta_{W_{i-1}} \label{eq:tautheta}\\
\gamma_i = \gamma_1\circ \theta_{W_{i-1}-\tau_1}
\label{eq:gammatheta}
\end{align}

From now on, we can write the strong Markov property, that we remind here
\begin{prop}[Strong Markov property]
\label{prop:strongMC}
For any Markov process $X$ regarding the filtration $(\mathcal{F}_n)$ with the initial distribution $\mu$ and any stopping time $T$ regarding the same filtration, if $\phi$ is $X$-measurable, then
\begin{equation}
 \mathbb{E}_{\mu}\left[\mathbb{I}_{\left\{T<\infty\right\}}\cdot\phi\circ\theta_T\middle|\mathcal{F}_{T}\right] = \mathbb{I}_{\left\{T<\infty\right\}}\cdot\mathbb{E}_{\left(X_T\right)}\left[\phi\right].
\end{equation}
\end{prop}
We apply this property to the relations \eqref{eq:tautheta} and \eqref{eq:gammatheta}, with $\psi = \mathbb{I}_{\left\{\tau_i=k\right\}}$ to obtain

\begin{align}
\mathbb{E}_{\mu_0}\left[\mathbb{I}_{\left\{W_{i-1}<\infty\right\}}\cdot\mathbb{I}_{\left\{\tau_i=k\right\}}\middle|\mathcal{F}_{W_{i-1}}\right] &= \mathbb{I}_{\left\{W_{i-1}<\infty\right\}}\cdot\mathbb{E}_{\left(H_{W_{i-1}}\right)}\left[\mathbb{I}_{\left\{\tau_1=k\right\}}\right] \\
&=\mathbb{P}_{X_{W_{i-1}}=Y_{W_{i-1}}}\left(\tau_1=k\right)\label{eq:taumarkov}
\end{align}
Yet 
$$
\mathbb{E}_{\mu_0}\left[\mathbb{I}_{\left\{W_{i-1}<\infty\right\}}\cdot\mathbb{I}_{\left\{\tau_i=k\right\}}\middle|\mathcal{F}_{W_{i-1}}\right] = \mathbb{P}\left(\tau_i=k,  W_{i-1}<\infty\middle| W_{i-1}\right)
$$
so
\begin{align*}
\mathbb{P}\left(\tau_i=k,  W_{i-1}<\infty\middle| W_{i-1}\right) &= \mathbb{P}\left(\tau_i=k\middle| W_{i-1}\right) \text{  a.s.}\\
&= \mathbb{P}_{H_{W_{i-1}}}\left(\tau_1=k\right) \text{  a.s.}
\end{align*}
Similarly, we can obtain
$$
\mathbb{P}\left(\gamma_i=k\middle| W_{i-1}\right) = \mathbb{P}_{H_{\tau_1}}(\gamma_1=k) \text{  a.s.}
$$
We still do the approximation that, for $\epsilon$ small enough, $\tau_i$ follow a geometric distribution of parameter $\epsilon$, insofar as the $\tau_i$ are now independent and identically distributed. As soon as we have $X_{W_{i-1}} = Y_{W_{i-1}}$, $\tau_i$ doesn't depend on the distribution of $(X_{W_{i-1}},Y_{W_{i-1}}$.\\
Regarding $\gamma_i$, we have to pay attention to the distribution of $H_{\tau_1}$. With the definition of $\tau_1$, we know that $X_{\tau_1}\neq Y_{\tau_1}$, but $X_{\tau_1-1}=Y_{\tau_1-1}$. We also know that the desynchronization that gives the time $\tau$ is due to the first time when the RMS doesn't follow the RDS transition, which happens with a probability close to $\epsilon$. So, on average, the two paths have been synchronized for approximately $\frac{1}{\epsilon}$ steps, following the mean MC behavior. We can so expect that, still on average, the distribution of $\left(X_{\tau_1-1},Y_{\tau_1-1}\right)$ is given by $(\pi)$ such that $\mathbb{P}\left(X_{\tau_1-1}=Y_{\tau_1-1}=i\right) = \pi(i)$. With this approximation, the distribution of $(X_{\tau_1},Y_{\tau_1})$ is given by $\mu_1$ (as defined in \eqref{eq:proba_desynch}), and so it is the same for any $(X_{\tau_i},Y_{\tau_i})$. The same argument also shows that the $\gamma_i$ are independent one from each other, as they only depend on their initial distribution $\mu_1$.\\

\subsection{Renewal Process results}

As we have seen, the couples $(\tau_i,T_i)$ are independent one from each other, with the same initial distribution. The process is so satisfying the properties of the renewal processes, as described in \cite{Taylor1998}, Chapter VII, and more precisely the paragraph (5.3.), for ``Renewal processes involving two components to each renewal interval'', also called ``Alternating renewal processes'' (more information and demonstration on \cite{random}). We have the main result
\begin{prop}[Asymptotic probability of synchronized time]
Let $p(t)$ be the probability that, at any step $t$, the process $H_t$ is synchronized (i.e. in a state of $S$). Then, the limit behavior of $p$ is
\begin{equation}
\lim_{t\to\infty} p(t) = \frac{\mathbb{E}_{\mu_1}\left[\tau_2\right]}{\mathbb{E}_{\mu_1}\left[T_2\right]} = \frac{\mathbb{E}_{\mu_1}\left[\tau_2\right]}{\mathbb{E}_{\mu_1}\left[\tau_2+\gamma_2\right]}.
\label{eq:synch_rate}
\end{equation}
\label{synchronized_prop}
\end{prop}
This means that, if one pick a time $t$ big enough, then the probability that the paths are synchronized at this time is given by the formula \eqref{eq:synch_rate}. So the expected rate of synchronized time is given by the same formula, which we can rewrite, with the approximation that the unsynchronization probability at each path is $\epsilon$, as
\begin{equation}
\frac{\mathbb{E}_{\mu_1}\left[\tau_2\right]}{\mathbb{E}_{\mu_1}\left[T_2\right]} = \frac{\frac{1}{\epsilon}}{\frac{1}{\epsilon}+\mathbb{E}_{\mu_1}\left[\gamma_2\right]} = \frac{1}{1+\epsilon \mathbb{E}_{\mu_1}\left[\gamma_2\right]}.
\label{eq:synch_rate2}
\end{equation}

\begin{figure}[h!]%
\includegraphics[width=\columnwidth]{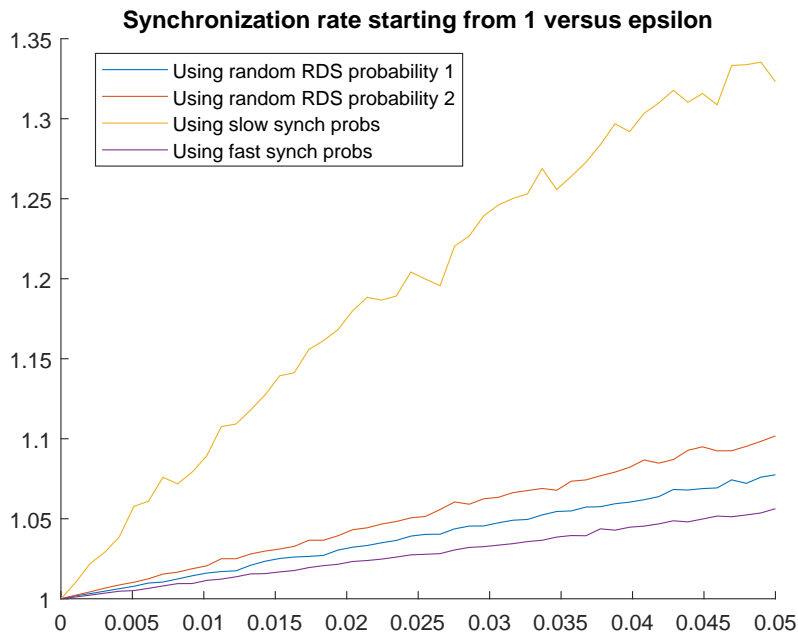}%
\caption{Unlike the three other ones, the yellow line, which represents the process with a long resynchronization time, is more chaotic, and does not behave as nicely as the other ones as a line.}%
\label{fig:synchro_rate2}%
\end{figure}

Let's draw some numerical realizations to evaluate how accurate this analysis is. We first study a two states system, so $s=2$, and as we have already explained, $Q$ is already determined by \eqref{eq:Q2}. The, for each probability distribution, we draw a realization of the comparison process RDS/RMS of $t=10^5$ steps, and count the number of steps where the RDS and the RMS are synchronized, i.e. the process is in state $S$, in order to calculate the synchronization rate (defined as the ratio of the synchronized time over the whole time). As the time $t$ is long, the ratios should follow the law given by \eqref{eq:synch_rate2} according to property \ref{synchronized_prop}. We make $\epsilon$ vary, and we plot the inverse of the ratio as a function to $\epsilon$, in order to see straight lines. In the figure \ref{fig:synchro_rate2}, we use different probability distribution: two are random, one has a small probability of drawing the synchronization matrices, and the last has a big probability of drawing them (in dimension two, there are only two synchronization matrices, $\left( \begin{array}{cc}
1 & 0\\
	1 & 0\\
\end{array}\right)$ and $\left( \begin{array}{cc}
0 & 1\\
	0 & 1\\
\end{array}\right)$, so, by changing their probability distribution, it is easy to chose a slow or fast resynchronization time $\gamma$).

\begin{figure}[h!]%
\includegraphics[width=\columnwidth]{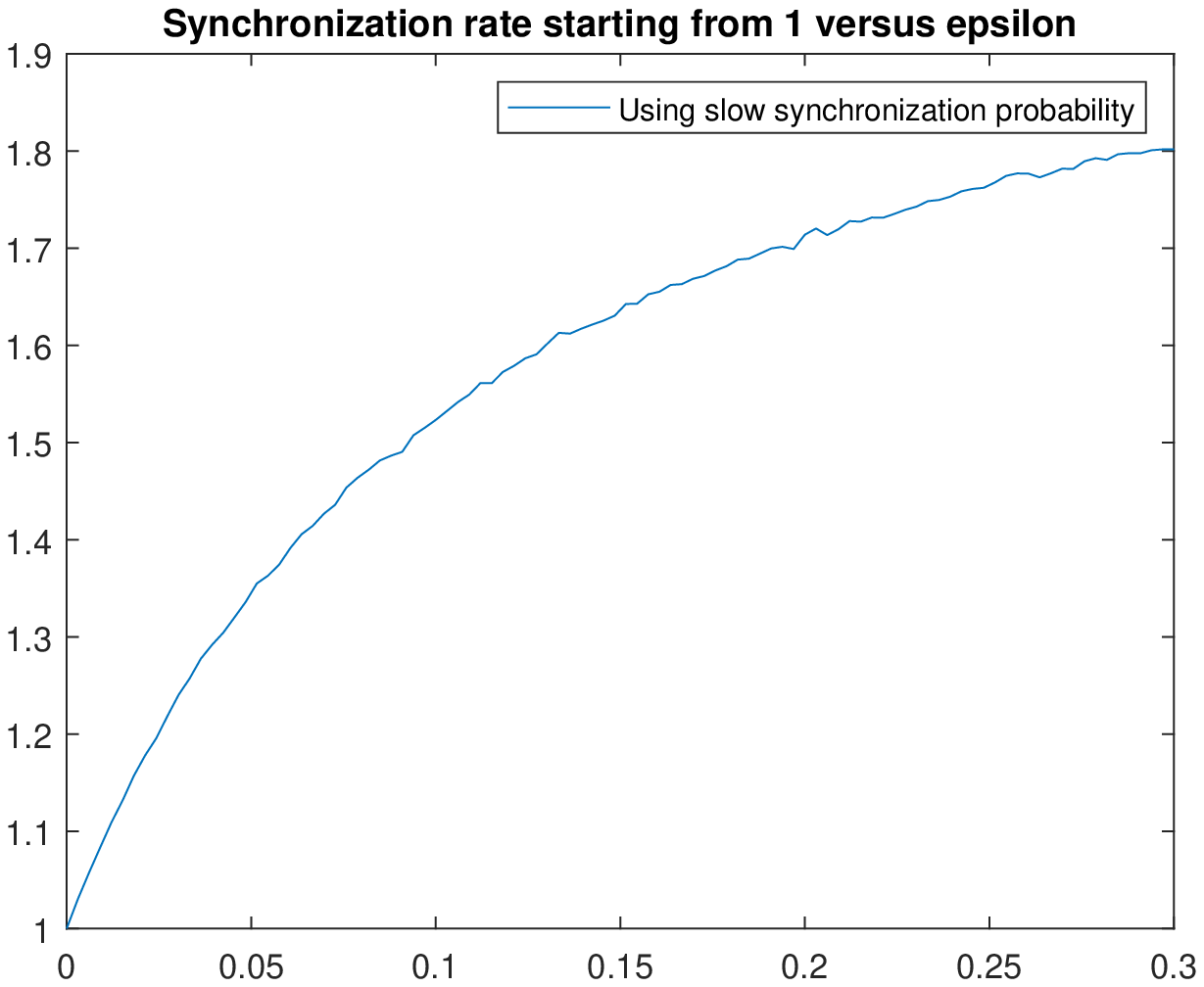}%
\caption{The line is smoother since the realization is $10^6$ steps long here, which reduces the variability of the results. Otherwise, the same parameters than the previous figure have been used.}%
\label{fig:synchro_rate_2_e6}%
\end{figure}

Except for the realization with a low probability of drawing the synchronizing matrices, all the plots behave very nicely as straight lines, and the slope is equal to the expected time of synchronization with the steady state distribution $\mu_1$ of the MC embedded in the RDS as initial distribution ($\mathbb{E}_{\mu_1}\left[\gamma_2\right]$). Let us spend more attention to the first path. The synchronization probability at each step was $\frac{1}{10}$, which leads to an expected synchronization time of $10$. In the meanwhile, if $\epsilon=0.05$, the RMS would not follow the RDS path with the rate of $1$ bifurcation over $20$ steps. The probability that the intrinsic noise prevents the synchronization is very small ($\frac{1}{200}$, drawing a synchronization matrix and not following the RDS at the same time), but that intrinsic noise could lead to a synchronization before the RDS, especially when the size $s$ of the state space is small. This phenomenon explains why this yellow line is so chaotic and under its slope very different from $10$ : the slope varies a lot with $\epsilon$, more than the processes with faster synchronization time. It is even more perceptible if we look at bigger values of $\epsilon$, as it is in figure \ref{fig:synchro_rate_2_e6}, where the straight line property is completely lost.

\begin{figure}[h!]%
\includegraphics[width=\columnwidth]{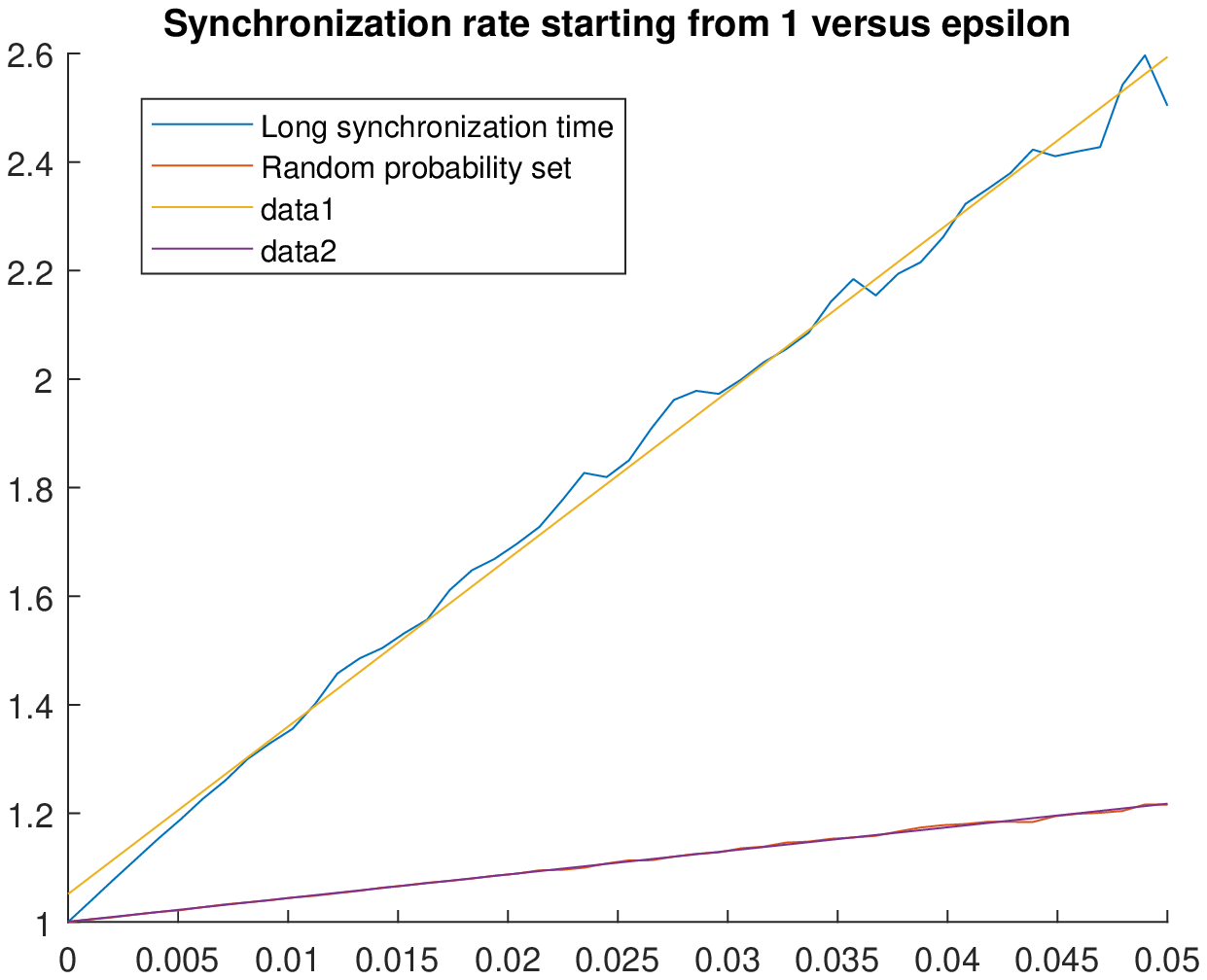}%
\caption{The chaotic line is produced by a distribution with a very low probability for partial synchronizating matrices and a higher one for invertible matrices. A random distribution over all the matrices produces the other one. The slope is particularly lower, so is the resynchronization time. The straight lines are the polynomial fits of degree one, and the realization is still $10^6$ steps long.}%
\label{fig:synchro_rate_3_e6}%
\end{figure}

Let now increase $s$ from $2$ to $3$. There are now $27$ different deterministic transition matrices, among which $3$ synchronize all the paths at the same time (matrices with one column of $1$) and $6$ that don't synchronize any path (the invertible ones, double stochastic matrices). So the $18$ left lead to ``partial synchronization'', they synchronize the paths if these later are coming from specific given states. And when $s$ gets bigger, the rate those matrices that conduct to partial synchronization keeps increasing, so that if we give a random probability distribution on the set of matrices, there is actually a bigger rate of ``synchronizing matrices'' that are weighted with a non negligible probability, and so the expected synchronization time for a random probability distribution is not very long. The figure \ref{fig:synchro_rate_3_e6} compares the inverse of the synchronization rate between a distribution with a slow synchronization and a random one. We still notice that the behavior is way smoother with a quick resynchronization, and the curve converges faster towards the straight line given by the formula.

Graphs + importance of $\epsilon$ in $\mathbb{E}_{\mu_1}\left[\gamma_2\right]$.
Dependance of $\epsilon$ of the slope should be less significant with higher $s$;

\section{Comments and further work}

\subsection{Critique}
\label{criticism}

To carry out the results of this work, we have had to use a certain number of hypothesis and approximations. Let us discuss a bit more about the relevance and accuracy of the two most important ones.\\
\newline
First, in order to use the nice results of renewal theory, we have considered that the random times involved in the process are independent and identically distributed. According to the precision and the consistency of the values (when we compare the slope of the straight lines in the synchronization rate versus $\epsilon$ and the analytical value calculated with the matrix $W$), this hypothesis does not seem to be nonsensical. Without getting into the regressed proof, one can easily understand that the $\tau_i$ are independent. $\tau_i$ is determined when the path of $Y$ diverges from the one of $X$, which happens with different probabilities according to the arriving state of the transition of $X$ - the probabilities which are the diagonal coefficients of $N$. It obviously does not depend on the last value of $\tau_{i-1}$ thus, but, if we consider that all of those coefficients of $N$ are non identical (we made the approximation that they are all equal to $1-\epsilon$), the value of $\tau_i$ depends on the distribution of the first state of the sequence - the first synchronized state after a desynchronized phase. This distribution may vary with $i$, and so the $(\tau_i)$ would be independent, but non identically distributed.\\
One can have the same reasoning for the $\gamma_i$ : the time duration of a synchronization phase does not depend of the length of the previous synchronization phase, but it definitely does depend on the initial state of that phase (which is the first non synchronized state after a synchronized phase). Our assumption is to consider that the synchronized phase is long enough to ensure that, ``on average'', the distribution of the last couple of synchronized states is given by the steady state distribution of the mean MC of the process (at least on the repetition of the following unsynchronization phases), and so we can use $\mu_1$ as the initial distribution of each sequence. There is a double assumption here actually. First, the mean MC has a steady state distribution and it converges towards it, second, the convergence towards this distribution is faster than the approximated duration of this phase $\frac{1}{\epsilon}$. As a MC on a finite state space, the steady state distribution necessarily exists, but the convergence would require some more properties (for example irreducibility and aperiodicity to fulfill Perron-Frobenius theorem requirements). A further study of strict RDS properties would gives us more information on the limit behavior of the path of a RDS regarding to its mean transition matrix, but it seems fair to assert this hypothesis.\\
\newline
Second, we want to discuss about the utility of our study. Everyone would agree that the main, persistent, most repeated hypothesis of this work is the fact that $\epsilon$ is ``small''. Let us remind the scale : small enough that the polynomial expansion of the exponential can be reduced at the $1$-st order in $\epsilon$, so we need $\epsilon$ to be negligible in comparison to $\epsilon^2$ (and obviously all of the upper orders). One may assume that $0.1$ is a good boundary, a pickier reader although might consider it under $0.05$ or $0.01$. It truly depends on the accuracy needed. However, considering this discussion, the equation \eqref{eq:synch_rate2} naturally shows that the synchronization rate converge to $1$ as $\epsilon$ gets smaller : the lighter the perturbation, the longer the synchronized time.\\
The real interest of that model lies in the relative values of $\epsilon$ and $\mathbb{E}_{\mu_1}\left[\gamma_1\right]$ such that their product has a not too close to zero value. For example, if $\mathbb{E}_{\mu_1}\left[\gamma_1\right] = \frac{1}{19}$, then the process will spend approximately $95$\% of the time synchronized, but if the product is equal to $2$, then it would only spend one third of the time synchronized.
 
\subsection{Related work: Mutual Information}

%Mutual information : definition + formula + related to further work

As you have already understood, the whole point of this paper was to understand how the introduction of this shape of intrinsic noise disturbs a RDS, and how would the new modified system would behave regarding the original one. There actually exists a more sophisticated mathematical tool that could be interesting to apply here : the mutual information. It can be seen as a way to measure the independence of two random variables as the information they share during their mutual realization.
\begin{mydef}{Mutual Information}
For two random variables $X$ and $Y$ taking their values into a discrete space $E$, then we define the mutual information $MI(X,Y)$ as
$$
MI(X,Y) = \sum_{x,y\in E}{\mathbb{P}\left(X=x, Y=y\right)ln\left(\frac{\mathbb{P}\left(X=x, Y=y\right)}{\mathbb{P}\left(X=x\right)\mathbb{P}\left(Y=y\right)}\right)}.
$$
The Mutual Information (MI) is a non-negative symmetric quantity, and is equal to $0$ only when the variables $X$ and $Y$ are independent.
\label{defmutualinf}
\end{mydef}
The name ``Mutual Information'' comes from the entropy theory, as we can write $MI(X,Y) = H(X) - H(X|Y) = H(X) + H(Y) - H(X,Y)$ when $H$ is an entropy function.

\begin{figure}[h!]%
\includegraphics[width=\columnwidth]{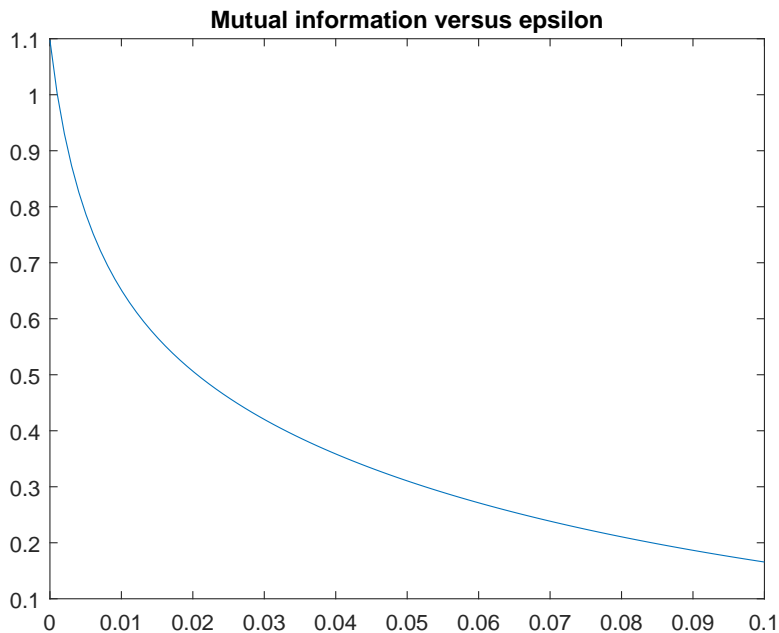}%
\caption{The MI is computed from the same slow resynchronization probability distribution over $\mathcal{D}$, in dimension $s=3$. For $\epsilon=0$, the paths are identical, and so the correlation between them is high. As $\epsilon$ gets higher, the RDS and the RMS are more and more independent, and so their MI decreases. Note that the decrease is more than exponential for very small $\epsilon$. The value of MI has been rescaled by the duration of the paths.}%
\label{mi_eps}%
\end{figure}

\begin{figure}[h!]%
\includegraphics[width=\columnwidth]{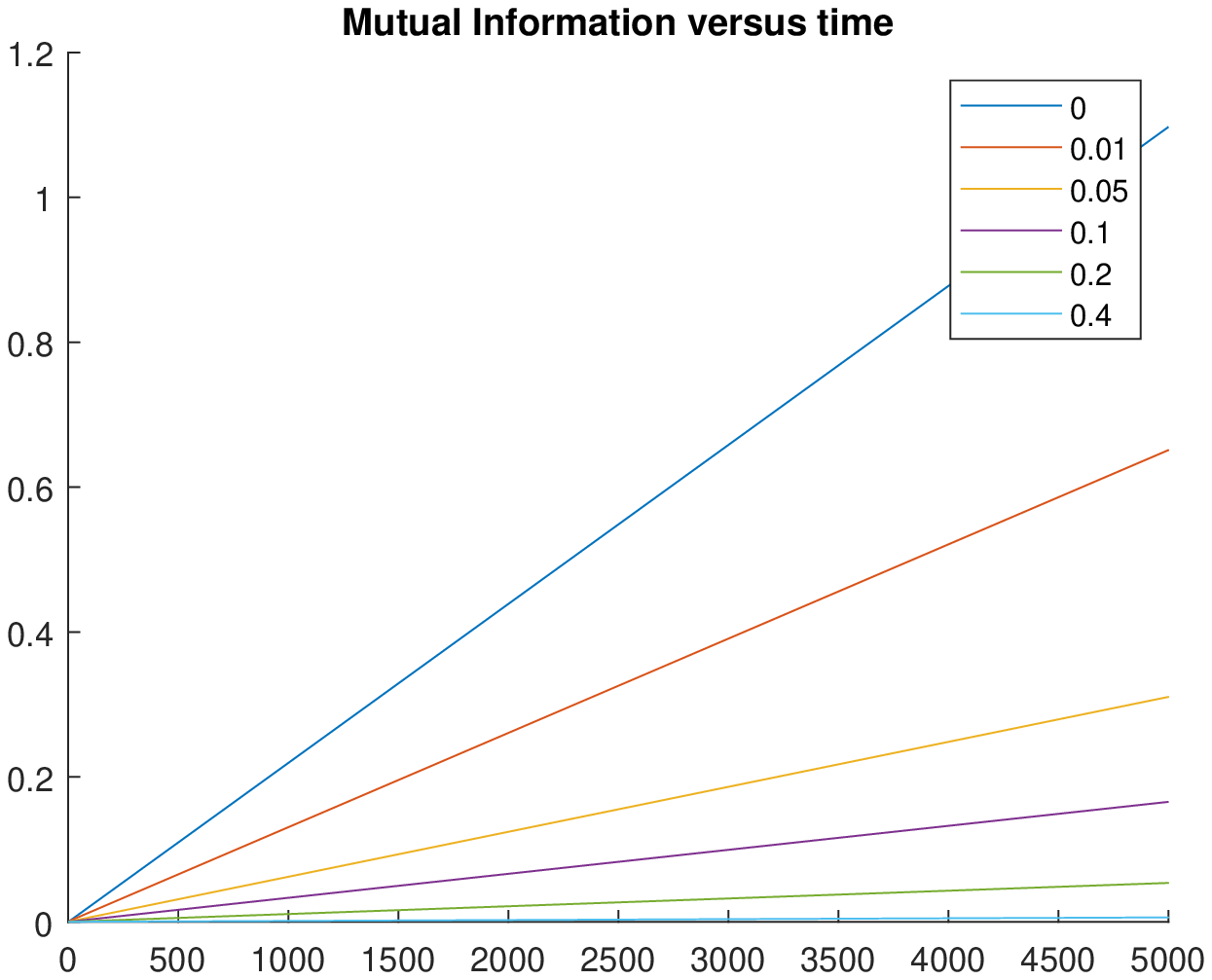}%
\caption{Same rescaling than the previous figure. Here the MI is plot versus the time, for different values of $\epsilon$. We notice the same behavior : the MI gets smaller as $\epsilon$ gets higher. More interestingly, the increase through the time is linear.}%
\label{mi_time}%
\end{figure}

In our situation, the comparative study of the RDS and the RMS are time dependent, so either we try to apply the MI definition on the current state of the process or on the previous states until time the current time, considering the sequence of the states as a random variable. The following application handles the second option. Let us call $\alpha_n$ and $\beta_n$ the sequences of states respectively from the processes $X$ and $Y$ until time $n$. We define the MI up to time $n$ starting from $(x_0,y_0)$ as
\begin{equation}
MI^{\left(n\right)}(x_0,y_0) = MI(\alpha_n,\beta_n).
\label{eq:mutualinf}
\end{equation}
with the definition of $MI$ given earlier and where the paths start respectively at states $x_0$ and $y_0$. We obviously use the mean transition probabilities given by $M$ and $MN$ to address the calculus of this quantity. By separating $MI^{\left(n\right)}$, a few lines of calculus lead to the formula
\begin{equation}
MI^{\left(n+1\right)}(x_0,y_0) = MI^{\left(n\right)}(x_0,y_0) + \sum_{k,l \in \mathcal{S}}{W_{(x_0,y_0)\rightarrow (k,l)}^n MI^{\left(1\right)}(k,l)}
\label{eq:recMI}
\end{equation}
using the double point motion mean transition matrix $W$ defined in the section \ref{transition_mat}. This formula is in practice very useful as it allows us to compute numerically the MI through the time. Indeed, one only need to initially calculate $W$ and $MI^{\left(1\right)}$ (stored as a matrix), and recursively apply the formula. The results are displayed on figure \ref{mi_eps} and \ref{mi_time}.

One can notice that the MI increases linearly with the time, so the crucial quantity that describes the MI would be the slope of the line. Analytically, as $W$ is a Markovian matrix, it might converge (under some conditions) towards a matrix $W^\infty$, and so the formula \eqref{eq:recMI} would describe the slope as $\sum_{k,l \in \mathcal{S}}{W^\infty_{(x_0,y_0)\rightarrow (k,l)} MI^{\left(1\right)}(k,l)}$.

Further work can be done to understand more deeply the link between the independence of the RDS/RMS and the parameters involved in the model (especially $\epsilon$ and the probability distribution $(q_i)$).

\section{Conclusion}
%Summarizing + further work : intrinsic noised induced synchro\\

Random dynamical systems are getting better and better known, as they embody a way to model extrinsic noise to a system of different entities. The point of this work was to study a way to introduce intrinsic noise in this particular structure while keeping the properties of the extrinsic noise.

Here, we have presented a way to do so, using a few parameters as the intensity of this noise $\epsilon$ or the perturbation matrix $Q$. We have come up with an adapted mathematical environment in order to have a regressed study of the phenomenon, and understand how it is influenced by its parameters. Despite the few analytical results, we have had a good understanding of this new tool, which opens door to a more subtle analysis of the relative weight of the intrinsic and extrinsic noises.

It needs now to be put into practice with some network model to understand towards which direction the work should be continued, however these random synchronization phases might have many applications in different fields. Up to now, the extrinsic was leading the synchronization phenomenon, while the intrinsic noise would cause the separation of the paths. But we can also think in a model where the extrinsic noise does not induce any synchronization (i.e. the RDS does not synchronize, only invertible matrices), and so some synchronization would happen because of the intrinsic noise, and persist during the time thanks to the extrinsic noise structure. That way, the choice of the parameters is crucial, and may lead to very different behaviors.

\bibliographystyle{plain}
\bibliography{C:/Users/Adrian/Documents/biblio}

\begin{thebibliography}{1}

\bibitem{Arnold1991}
Ludwig Arnold and Hans Crauel.
\newblock Random dynamical systems.
\newblock In Ludwig Arnold, Hans Crauel, and Jean-Pierre Eckmann, editors, {\em
  Lyapunov Exponents}, pages 1--22, Berlin, Heidelberg, 1991. Springer Berlin
  Heidelberg.

\bibitem{mcinre}
Robert Cogburn.
\newblock Markov chains in random environments: The case of markovian
  environments.
\newblock {\em The Annals of Probability}, 8(5):908--916, 1980.

\bibitem{MDP1970}
C.~Derman.
\newblock {\em Finite State Markov Decision Processes}.
\newblock Academic Press, 1970.

\bibitem{random}
Kyle Siegrist.
\newblock Probability, mathematical statistics, stochastic processes.

\bibitem{takahashi1969}
Yukio Takahashi.
\newblock Markov chains with random transition matrices.
\newblock {\em Kodai Math. Sem. Rep.}, 21(4):426--447, 1969.

\bibitem{Taylor1998}
H.~M. Taylor and S.~Karlin.
\newblock {\em An Introduction to Stochastic Modeling}.
\newblock Academic Press, 1998.

\bibitem{X.-F.Ye2016}
Felix X.-F.~Ye, Yue Wang, and Hong Qian.
\newblock {S}tochastic dynamics: {M}arkov chains and random transformations.
\newblock {\em Discrete and Continuous Dynamical Systems - Series B},
  21:2337--2361, 2016.

\end{thebibliography}

\end{document}